\documentclass[reqno,12pt]{amsart}
\usepackage{latexsym,amsmath,amsthm,amssymb}

\newtheorem{theorem}{Theorem}[section]
\newtheorem{lemma}[theorem]{Lemma}
\newtheorem{corollary}[theorem]{Corollary}%[section]
\newtheorem{proposition}[theorem]{Proposition}%[section]
\theoremstyle{remark}
\newtheorem{example}[theorem]{Example}
\newtheorem{remark}[theorem]{Remark}%[section]
\theoremstyle{definition}

\numberwithin{equation}{section}

\begin{document}

\title{Mass transport generated by a flow of Gauss maps}
\author{Vladimir I. Bogachev, Alexander V. Kolesnikov}

\begin{abstract}
Let $A \subset \mathbb{R}^d$, $d\ge 2$, be a compact convex set   and let
$\mu = \varrho_0\, dx$ be a probability
measure on $A$ equivalent to the restriction of Lebesgue measure.
Let $\nu = \varrho_1\, dx$ be a probability measure
on $B_r := \{x\colon\, |x| \le r\}$ equivalent to the restriction
of Lebesgue measure.
We prove that there exists a mapping $T$ such that  $\nu = \mu \circ T^{-1}$
and $T = \varphi \cdot {\rm n}$, where $\varphi\colon\, A \to [0,r]$
is a continuous potential with convex sub-level sets
and ${\rm n}$ is the Gauss map of the corresponding level sets
of~$\varphi$. Moreover, $T$ is invertible and essentially
unique. Our proof employs
the optimal transportation techniques.
We show that
in the case of smooth $\varphi$ the level sets of $\varphi$ are
governed by the Gauss curvature flow
$\dot{x}(s) = -s^{d-1} \ \frac{\varrho_1(s  {\rm n})}{\varrho_0(x)}  K(x)
\cdot {\rm n}(x)$, where $K$ is the Gauss curvature.
As a by-product one can reprove the existence of
weak solutions to the classical Gauss curvature
flow starting from a convex hypersurface.
\end{abstract}

\maketitle



\section{Introduction}

The goal of this paper is to introduce a new class of
transformations of measures on $\mathbb{R}^d$ which
(heuristically) have the form
$T=\varphi \cdot \nabla\varphi/|\nabla\varphi|$
with some function~$\varphi$.
Our work is motivated by  two intensively developing areas:
optimal transportation and curvature flows, and establishes an
interesting link between these areas.
Optimal transportation can be described as a problem
of optimization of a certain functional associated
with a pair of measures.
The quadratic transportation cost   $W_2^2(\mu, \nu) $  between two
probability measures $\mu, \nu$ on $\mathbb{R}^d$ is
 defined as the minimum of the following functional (the
Kantorovich functional):
\begin{equation}
\label{MKpr}
m \mapsto \int_{\mathbb{R}^d\times \mathbb{R}^d}
|x_1 - x_2|^2 \,d m(x_1,x_2), \quad m\in \mathcal{P}(\mu,\nu),
\end{equation}
where $\mathcal{P}(\mu,\nu)$ is the set of all probability measures
on $\mathbb{R}^d\times\mathbb{R}^d$ with the marginals $\mu$ and~$\nu$;
here $|v|$ denotes the Euclidean norm of $v\in \mathbb{R}^d$. The problem of
minimizing (\ref{MKpr}) is called the mass transportation problem.
This formulation is due to Kantorovich~\cite{Kant}.
A~detailed discussion of the mass
transportation problem in this setting
can be found in \cite{RR}.
In many cases there exists a mapping
$T\colon\, \mathbb{R}^d \to \mathbb{R}^d$, called
the optimal transport between $\mu$ and $\nu$ (or a solution to the Monge problem),
such that $\nu = \mu \circ T^{-1}$ and
$$
W^2_2(\mu, \nu) =\int_{\mathbb{R}^d} |x -T(x)|^2 \, \mu(dx).
$$
The minimization of the latter integral in the class
of measurable mappings $T$ such that $\mu\circ T^{-1}=\nu$ is called
the Monge problem. If $T$ is a solution to the Monge problem,
then the image of $\mu$ under the mapping $x\mapsto (x,T(x))$
to $\mathbb{R}^d\times \mathbb{R}^d$ minimizes the Kantorovich
functional.
However, it may happen that the Monge problem has no solution,
while the Kantorovich problem is always solvable.
It is worth mentioning that
the first rigorous results related to existence of optimal
mappings were obtained in the classical work of A.D.~Alexandroff
\cite{Al} on convex surfaces with prescribed curvature!
If $\mu$ and $\nu$ are absolutely continuous, then, as show by
Brenier  \cite{Bren} and McCann \cite{McCann},
 there exists an optimal transportation $T$ which takes
 $\mu$ to~$\nu$. Moreover,
this mapping is $\mu$-unique  and has the form $T =
\nabla W$, where $W$ is convex. Under broad assumptions,  $W$  solves the
following nonlinear PDE (the  Monge--Amp{\`e}re equation):
$$
\varrho_\nu(\nabla W)
\det D^2_a W = \varrho_\mu,
$$
where
$\varrho_\mu$ and $\varrho_\nu$ are densities of $\mu$ and $\nu$
and
$D^2_a W$ is the absolutely continuous part of the distributional derivative of $D^2 W$.
At present the optimal transportation theory  attracts attention of researchers from
the most diverse fields, including probability, partial differential
equations, geometry, and infinite-dimensional analysis (see
Villani's book \cite{Vill} and papers
\cite{Ambrosio} and~\cite{Sturm}).

The study of curvature flows is a very popular subject in geometry.
The theory of Ricci flows attracted particular interest
after the famous works of G.~Perelman  on the Poincar{\'e} conjecture.
The theory of geometrical flows began, however, with flows of embedded manifolds.
Let $F_0 \colon\, M^{d-1} \to \mathbb{R}^d$ be a smooth embedding of a smooth compact
Riemannian manifold $M^{d-1}$ (without boundary).
Denote by   $A$ the enclosed body:  $\partial A = F_0(M)$.
We say that an evolution
$$
F(\,\cdot\,,\,\cdot\,) \colon\, M \times [0,T) \to \mathbb{R}^d
$$
is a geometrical flow if  $F_0 = F(\cdot,0)$ and
\begin{equation}
\label{meancurv}
\frac{\partial}{\partial t} F(x, t)
= - g(F(x,t))\cdot  {\rm n}(F(x,t)),
\end{equation}
where $g\colon\, M \to \mathbb{R}$ is
some curvature function  and  ${\rm n}$ is the outer unit normal vector.
If $g=H$ is the mean curvature, then $F$ is called  the mean curvature flow.
If $g=K$ is the Gauss curvature, then $F$ is called the Gauss curvature flow.
The simplest case of the  curvature flow is given by the following planar flow:
$$
\frac{\partial x(t,s)}{\partial t}  = - k(x) \cdot {\rm n}(x),
\ x \in \mathbb{R}^2, \ M=S^1.
$$
Consider the flow of closed curves $t \mapsto x(t,\cdot)$.
Under this flow the enclosed volume
decreases with the constant speed $-2\pi$.
In addition,
any non-convex curve becomes convex in finite time and
then remains convex.
Finally,  any curve
shrinks to a point in finite time; the shape of any curve
  becomes more and more rotund
  (see \cite{GaHa} and \cite{Gray}).
Any multi-dimensional
mean curvature flow or Gauss curvature flow
starting from a convex surface preserves convexity and
shrinks the surface to a point (see \cite{Huis} and
\cite{Tso}).
In \cite{EvSp} and \cite{CGG},
equation (\ref{meancurv})  in the case of the
mean curvature was investigated from the PDE's point of view.
It turns out that
the surfaces driven by (\ref{meancurv})
can be obtained as level sets of a function $u(t,x)$ which
satisfies  a  nonlinear degenerate second order parabolic equation of
the Monge--Amp{\`e}re type.
A solution of this equation is in general understood in some weak sense (viscosity solutions).
For the PDE approach and viscosity solutions,
see the recent book \cite{Giga}.
Concerning Gauss flows, see \cite{BenAnd}.

In this paper, we establish the existence and uniqueness
of a special  measure
transportation mapping
between two probability measures $\mu$ and~$\nu$.
It has the following heuristic expression:
$$
T = \varphi \frac{\nabla \varphi}{|\nabla \varphi|}.
$$
The potential $\varphi$ has convex sub-level sets.
Note that this transportation mapping may not be a gradient.
Nevertheless, $T$ can be obtained as a degenerate limit of some transportation mappings
which are constructed by means of the optimal transportation techniques.
The limiting potential
satisfies a degenerate Monge--Amp{\`e}re
equation (see the proof of the main theorem).
In addition, we show that the resulting limit is naturally connected with the
Gauss flow. The level sets of  the potential $\varphi$ can be associated to
a special Gauss flow
associated to the measures $\mu$ and $\nu$ according to
$$
\dot{x}(s) = -s^{d-1} \ \frac{\varrho_1(s  {\rm n})}{\varrho_0(x)}
K(x)\cdot {\rm n}(x).
$$
In the case $\varrho_1(x) = \frac{C_{d,r}}{|x|^{d-1}}$ and
$\varrho_0(x) = \frac{1}{\mathcal{H}^d(A)}$
we obtain a weak solution to
$$
\dot{x}(s) = - c \ K(x)\cdot {\rm n}(x),
$$
which is the classical Gauss flow starting from some initial convex hypersurface.

Finally, we note that in
\cite{LottV},
\cite{McT},
\cite{Top2}
the reader can find
other interesting links between
 mass transportation  and  geometrical flows (in particular, the Ricci flows).
Some analogs of the presented results in the case of a manifold
will be considered in our forthcoming joint paper with F.-Y.~Wang.

\section{Main result}

Throughout we assume that $d\ge 2$ and denote by
$\mathcal{H}^{n}$ the
$n$-dimensional Hausdorff measure. For Lebesgue measure we often use
another common notation~$dx$. Let ${\rm Int}\, A$ denote the
interior of a set~$A$.

Recall that, given a compact smooth orientable $(d-1)$-dimensional surface
$M$ in $\mathbb{R}^d$, one has the Gauss map ${\rm n}\colon\, M \to S^{d-1}$, where
${\rm n}(x)$ is the global unit outer normal vector field. Let
$D{\rm n}\colon\, TM_x \to TS^{d-1}_{{\rm n}(x)}$ be
the differential of~$\rm{n}$. Choose an orthonormal basis
 $\{e_2, \ldots, e_d\} \subset TM_x$. Then the matrix $D \rm{n}$
 can be written as $\langle \partial_{e_i} {\rm{n}}, e_j \rangle$, where
 $\partial_{e_i} \rm{n}$ are the usual
 partial derivatives of $\rm{n}$.
The determinant of $D {\rm n}(x)$ is called the Gauss curvature
and is denoted throughout   by $K(x)$.

Below we deal with the case where $M$ is a
  surface (possibly, non-smooth) of the form
$M = \partial V$, where $V$ is a convex compact set.
In this case the normal ${\rm n}(x)$ is well-defined almost
everywhere on~$M$. More precisely,
for an arbitrary point $x \in M$, let us set
$$
N_{M,x} := \{\eta \in S^{d-1} ; \ \forall z \in V, \ \langle \eta, z-x \rangle \le 0  \}.
$$
If $N_{M,x}$ contains a single element ${\rm n}(x)$, then
${\rm n}(x)$ is the unit normal in the usual sense.
We shall use the fact that one has ${\mathcal H}^{d-1}(S)=0$, where
$$
S= \{x\colon\, N_{M,x} \ \hbox{contains more than one element}\}.
$$
Hence the Gauss map ${\rm n}(x)$ is well-defined
$\mathcal{H}^{d-1}$-almost everywhere on~$M$.
Moreover,  one can show that
$K(x)$ is well-defined
$\mathcal{H}^{d-1}$-almost everywhere on $M$
(but this fact is not used below).

 We shall consider the following Hausdorff distance between
 nonempty compact sets:
$$
{\rm{dist}}(B_1,B_2) = \max \bigl(\sup_{x \in B_1} {\rm dist}(x,B_2), \
 \sup_{x \in B_2} {\rm dist}(x,B_1) \bigr) .
$$

\begin{theorem}
\label{maintheorem}
Let $A \subset \mathbb{R}^d$ be a compact convex set   and let
$\mu = \varrho_0\, dx$ be a probability
measure on $A$ equivalent to the restriction of Lebesgue measure.
Let $\nu = \varrho_1\, dx$ be a probability measure on
$B_r = \{x\colon\, |x| \le r\}$ equivalent to the restriction
of Lebesgue measure.
Then, there exist a Borel mapping $T\colon\, A \to B_r$
and a continuous function $\varphi\colon\, A \to [0,r]$ with convex
sub-level sets
$A_s = \{\varphi \le s\}$
 such that $\nu = \mu \circ T^{-1}$
and
$$
T = \varphi \cdot {\rm n}
\quad
\hbox{$\mathcal{H}^d$-almost everywhere,}
$$
  where
${\rm n} = {\rm n}(x)$ is a
unit outer normal vector to the level set
$\{y\colon\, \varphi(y)=\varphi(x)\}$ at the point~$x$.

If $\varphi$ is smooth, the level sets of $\varphi$
are moving according to
the following Gauss curvature flow equation:
\begin{equation}
\label{GaussFlow}
\dot{x}(s) = -s^{d-1} \ \frac{\varrho_1(s {\rm n})}{\varrho_0(x)}  K(x)\cdot {\rm n}(x)
\end{equation}
where $x(s) \in \partial A_{r-s}$,  $0 \le s \le r$,
$x(0) \in \partial A$ is any initial point.
\end{theorem}

To prove this theorem we develop an approach based on the
optimal transportation techniques.
For every $t\ge 0$, we consider  a mapping
$T_t$ that takes $\mu$ to $\nu$
and maximizes  the functional
\begin{equation}
\label{MKt}
F \mapsto \int  \langle x, F(x)\rangle |F(x)|^t\, \mu(dx)
\end{equation}
in the class of mappings $F$ with $\mu\circ F^{-1}=\nu$.
Equivalently, it minimizes the functional
$$
F \mapsto \int  \bigl| x - F(x)|F(x)|^t \bigr|^2\, \mu(dx)
$$
in the class of mappings $F$ with $\mu\circ F^{-1}=\nu$.
For $t=0$ (\ref{MKt}) becomes the classical Monge--Kantorovich problem.
For $t \ne 0$  standard arguments from the Monge--Kantorovich theory
show that the set
$$
\Bigl\{ \bigl(x, T_t(x)|T_t(x)|^t \bigr), \ x \in A \Bigr\}
$$
is cyclically monotone, hence belongs to the graph of the
gradient of some convex function $W_t$
(see  \cite[Chapter~2]{Vill}).
This can be shown, for instance, by a cyclical permutation of small balls
(see \cite{Vill}). More formally, this can be obtained
by variation of the corresponding  Lagrange functional (see \cite{Evans3}).

If the reader does not want to be concerned with the cyclical monotonicity
or calculus of variations, we  note
that $\nabla W_t$ is just the optimal transportation of
$\mu$ to $\nu \circ S^{-1}_t$, where $S_t(x) = x |x|^t$.
This can be taken for a definition of $W_t$.

One has the following relations:
$$
T_t = \frac{\nabla W_t}{|\nabla W_t|^{\frac{t}{1+t}}}, \
\nabla W_t(x) = T_t(x)  |T_t(x)|^t.
$$
Clearly, $|T_t(x)|\le r$ since $T_t$ transforms $\mu$ into $\nu$.

Throughout the paper we choose $W_t$ in such a way that $\min_{x \in A} W_t(x) =0$.
Define  a new potential function $\varphi_t$ by
$$
W_t = \frac{1}{t+2} \varphi^{t+2}_t.
$$
One has
$$
T_t = \varphi_t \frac{\nabla \varphi_t}{|\nabla \varphi_t|^{\frac{t}{t+1}}}.
$$
We show below that the limits
$$
\lim_{t\to \infty}\varphi_t = \varphi, \
\lim_{t\to \infty}T_t = T
$$
exist almost everywhere (for a suitable sequence $t_n\to\infty$)
and then we prove that
$T$ is the desired mapping.

\begin{lemma}
\label{gradient-estimate}
One has
$$
\varphi_t \le
(2+t)^{\frac{1}{2+t}} \bigl({\rm diam}(A)\bigr)^{\frac{1}{2+t}}
 r^{\frac{1+t}{2+t}},
$$
$$
\int_{A} |\nabla \varphi_t(x)|\, dx   \le
\int_{\partial A} \varphi_t\, d\mathcal{H}^{d-1}
 \le
(2+t)^{\frac{1}{2+t}} \bigl({\rm diam}(A)\bigr)^{\frac{1}{2+t}}
 r^{\frac{1+t}{2+t}} \ \mathcal{H}^{d-1}(\partial A).
$$
\end{lemma}
\begin{proof}
By the convexity of $W_t$ we have
$$
W_t(x) - W_t(y) \le \langle x-y, \nabla W_t(x)\rangle.
$$
Choosing $y_0$ in such a way that $W_t(y_0) =0$, we find
$W_t(x) \le {\rm diam}(A) |\nabla W_t(x)|$ for every $x \in A$.
Since $W_t = \frac{1}{2+t} \varphi^{t+2}_t$, we obtain
$$
\varphi_t \le  (t+2) {\rm diam}(A) |\nabla \varphi_t|.
$$
Let $\alpha = \frac{1+t}{2+t}$.
Note that $\frac{1-\alpha}{\alpha} = \frac{1}{1+t}$, hence
$\varphi_t |\nabla \varphi_t|^{\frac{1-\alpha}{\alpha}}
= |T_t| \le r$. Therefore, one has
\begin{align*}
 \varphi_t = \varphi_t^{1-\alpha} \varphi_t^{\alpha}
&
\le \bigl((2+t) {\rm diam}(A)\bigr)^{1-\alpha}
  |\nabla \varphi_t|^{1-\alpha} \varphi_t^{\alpha}
\\&
\le
\bigl((2+t) {\rm diam}(A)\bigr)^{1-\alpha}
 \Bigl[ \varphi_t |\nabla \varphi_t|^{\frac{1-\alpha}{\alpha}} \Bigr]^{\alpha}
 \\&
 \le
 \bigl((2+t) {\rm diam}(A)\bigr)^{1-\alpha}
 r^{\alpha}
 =
 \bigl((2+t) {\rm diam}(A)\bigr)^{\frac{1}{2+t}}
 r^{\frac{1+t}{2+t}} .
\end{align*}
By using the convexity of $W_t$ again, we get
$$
0 \le {\rm div} \Bigl(\frac{\nabla W_t}{|\nabla W_t|}\Bigr)
= {\rm div} \Bigl(\frac{\nabla \varphi_t}{|\nabla \varphi_t|}\Bigr),
$$
where under ${\rm div} \Bigl(\frac{\nabla W_t}{|\nabla W_t|}\Bigr)$ we understand
the distributional derivative of the
vector field $\frac{\nabla W_t}{|\nabla W_t|}$.
Integrating with respect to $\varphi_t dx$ we obtain
$$
0 \le \int_{A} {\rm div}
\Bigl( \frac{\nabla \varphi_t}{|\nabla \varphi_t|} \Bigr)\varphi_t\, dx
= - \int_{A} |\nabla \varphi_t|\, dx
+ \int_{\partial A} \varphi_t \Bigl\langle n_A, \frac{\nabla \varphi_t}
{|\nabla \varphi_t|}\Bigr\rangle\, d\mathcal{H}^{d-1}.
$$
Hence
$$
\int_{A} |\nabla \varphi_t|\, dx
\le
\int_{\partial A} \varphi_t\,  d\mathcal{H}^{d-1}.
$$
Applying the above uniform estimate for $\varphi_t$ we complete the
proof.
In fact, one could do these calculations in the case of smooth densities,
where $\varphi_t$ has a better regularity, and then approximate
our densities by smooth ones (the corresponding optimal transports
converge to~$\nabla W_t$).
\end{proof}

\begin{corollary}
\label{potential-conv}
There exists a sequence $\{t_n\} \to \infty$ such that $\{\varphi_{t_n}\}$
converges almost everywhere to a finite function $\varphi$.
\end{corollary}
\begin{proof}
By Lemma \ref{gradient-estimate},  every sequence $\{\varphi_{t_n}\}$
is bounded in $W^{1,1}(A)$. By the compactness of the embedding
$W^{1,1}(B) \subset L^1(B)$ for any ball $B \subset A$ and the
diagonal argument (in fact, since in the
present situation $A$ is convex,
the embedding of the whole space $W^{1,1}(A)$ is compact),
we obtain the claim.
\end{proof}

\begin{lemma}
\label{powergrad}
There exists a sequence $\{t_n\} \to \infty$ such that
$$
\lim\limits_{t_n \to \infty}  |\nabla \varphi_n|^{\frac{1}{1+t_n}} =1
$$
almost everywhere, where $\varphi_n :=\varphi_{t_n}$.
\end{lemma}
\begin{proof}
By Lemma \ref{gradient-estimate} one has
$\varphi_t \le  (t+2) {\rm diam(A)} |\nabla \varphi_t|$.
Hence
$$
C_t |T_t|^{\frac{1}{2+t}}
\le
|\nabla \varphi_t|^{\frac{1}{1+t}},
$$
where
$C^{-1}_t = (2+t)^{\frac{1}{2+t}} \bigl({\rm diam}(A)\bigr)^{\frac{1}{2+t}}$.
Changing variables one gets the following estimate for any $\delta>0$:
$$
\mu \Bigl( C_{t_n}^{-1} |\nabla \varphi_n|^{\frac{1}{1+t_n}}
\le 1-\delta \Bigr)
 \le
\mu \Bigl(  |T_{t_n}|^{\frac{1}{2+t_n}}\le 1-\delta \Bigr)
 =
\nu \Bigl( |x|^{\frac{1}{2+t_n}}
\le 1-\delta\Bigr).
$$
Hence
$$
\mu \Bigl( 1- C_{t_n}^{-1} |\nabla \varphi_n|^{\frac{1}{1+t_n}}
\ge \delta \Bigr) \to 0.
$$
This implies that
$\bigl(1- C_{t_n}^{-1} |\nabla \varphi_n|^{\frac{1}{1+t_n}}\bigr)^{+}$
tends to zero in $\mu$-measure as $t_n \to \infty$.
Passing to an almost everywhere  convergent subsequence
one can assume additionally that
\begin{equation}
\label{inflim}
\underline{\lim}_{t_n \to \infty}|\nabla \varphi_n|^{\frac{1}{1+t_n}}
\ge 1
\end{equation}
almost everywhere.
Since $\sup_t \|\nabla\varphi_t\|_{L^1(dx)} < \infty$
by Lemma \ref{gradient-estimate},
we see  that the sequence $\{|\nabla \varphi_n|^{\frac{1}{1+t_n}}\}$
is bounded in $L^p(A)$ for any $p<\infty$.
Moreover, by H{\"o}lder's inequality
$$
\overline{\lim}_{t_n \to \infty} \int_{A}
|\nabla \varphi_n|^{\frac{p}{1+t_n}} dx \le \mathcal{H}^d(A).
$$
Hence, choosing an $L^p(A)$-weakly convergent subsequence
$|\nabla \varphi_{n_m}|^{\frac{1}{1+t_{n_m}}} \to f$, one has
$$
\int_{A} f\, dx \le \mathcal{H}^d(A).
$$
On the other hand, (\ref{inflim}) and Fatou's lemma show
that $f\ge 1$ a.e., which yields
$$
\lim\limits_{t_{n_m} \to \infty} \int_{A}
|\nabla \varphi_{n_m}|^{\frac{p}{1+t_{n_m}}}\, dx =1.
$$
Hence
$|\nabla \varphi_{n_m}|^{\frac{1}{1+t_{n_m}}} \to 1$
in the norm of any $L^{p}(A)$, $p<\infty$.
Extracting again an almost everywhere convergent subsequence we get the claim.
\end{proof}

In what follows we set
$\varphi_n :=\varphi_{t_n}$ and assume that $\varphi_n \to \varphi$
and $|\nabla \varphi_n|^{\frac{1}{1+t_n}} \to 1$
almost everywhere.

\begin{lemma}
\label{set-conv}
Let $C_n \subset B_r$  be convex sets  such that
$I_{C_n} \to I_C$
almost everywhere. If $C$ is of positive measure, then
$
{\rm dist}(\partial C_n, \partial C) \to 0.
$
\end{lemma}
\begin{proof}
The set $C$ can be taken convex by letting
$I_C: = \underline{\lim}_n I_{C_n}$. We may assume that $r=1$.
It is known and readily verified by induction that every
 convex set $U\subset B_r$ with $\mathcal{H}^d(U)\ge \delta$
 contains a ball of volume at least $\kappa_1(d)\delta$,
 where $\kappa_1(d)$ depends only on~$d$.
Let $B$ be a ball of radius $\varepsilon>0$
centered at some point $x_0\in \partial U$.
Then
$$
\mathcal{H}^d(U\cap B)\ge \kappa_2(d)\varepsilon^{d}\delta,
$$
where $\kappa_2(d)$ depends only on~$d$.
It follows that, whenever $\mathcal{H}^d(C_n)\ge \mathcal{H}^d(C)/2$, one has
$\|I_C-I_{C_n}\|_{L^1}=\mathcal{H}^d(C\bigtriangleup C_n)\ge
2^{-1}\kappa_2(d)\mathcal{H}^d(C) {\rm dist}(\partial C_n, \partial C)^d$.
\end{proof}

Note that due to convexity one has
${\rm dist}(\partial C_n, \partial C)={\rm dist}(C_n,C)$.

\begin{lemma}
\label{uniform-conv}
The sequence of potentials
 $\varphi_{t_n}$ converges  to $\varphi$
 uniformly on~$A$. In particular,   $\varphi$
 is continuous and has  convex sub-level sets $A_s = \{y\colon\, \varphi(y) \le s \}$.
\end{lemma}
\begin{proof}
Clearly, it is sufficient to prove the claim for
a subsequence. As noted above,
one can assume, in addition, that
$|\nabla \varphi_n|^{\frac{1}{1+t_n}} \to 1$ almost everywhere.
Let us redefine $\varphi$ as follows: $\varphi: = \underline{\lim}_n  \varphi_{n}$.
Then  the sub-level sets $A_{s}$ of $\varphi$ are convex since the
corresponding sub-level sets
$
A_{s,n} = \{x\colon\, \varphi_n(x) \le s\}
$
 of $\varphi_{n}$ are convex.
Since $|T_n| = \varphi_n |\nabla \varphi_{n}|^{\frac{1}{1+t_n}}$, we have shown that
 $|T_n| \to \varphi$ almost everywhere. Then the equality
$\nu = \mu \circ T^{-1}_n$ yields that the
image  of $\mu$ under the mapping
$
x \mapsto \varphi(x) \in \mathbb{R}^{+},
$
denoted by $\mu_{\varphi} \in \mathcal{P}(\mathbb{R}^+)$,
coincides with $\nu_{|x|} \in \mathcal{P}(\mathbb{R}^+)$, where $\nu_{|x|}$
is the image of $\nu$
under the mapping $x \mapsto |x|$.
Due to our assumptions on~$\nu$,
this implies that
$\mu_{\varphi}$  has a strictly increasing continuous distribution function, i.e.,

1) $\mu(A_{s_1}) < \mu(A_{s_2})$ whenever $s_1 < s_2$,

2)
$\mu(\{\varphi=t\}) =0$ for all $t \in [0,r]$.

Note that  2) implies that $I_{A_{s,n}} \to I_{A_s}$
 almost everywhere for each~$s>0$.
By Lemma \ref{set-conv} we have
$$
{\rm dist} (\partial A_{s,n}, \partial {A_s}) \to 0,\quad s>0.
$$
Now, given $\varepsilon>0$, we divide $[0,r]$ by points
$s_1,\ldots,s_N$ with $|s_{i+1}-s_i|<\varepsilon$ and take
$\delta = \max_{i \le N} {\rm dist} (\partial A_{s_i}, \partial A_{s_{i+1}})$.
There exists $M$ such that
$$
{\rm dist} (\partial A_{s_i,n},
 \partial A_{s_i}) < \delta/2
 $$
 for every $i=1,\ldots,N$ and every $n>M$.
 This implies that $\sup_{x \in A}|\varphi_n(x) - \varphi(x)| \le 2\varepsilon$ for
 all $n \ge M$. Hence $\varphi_n \to \varphi$ uniformly.
 Since $\varphi_n$ are continuous as powers of convex functions,
 $\varphi$ is continuous as well. The proof is complete.
\end{proof}

\begin{lemma}
\label{normalcone}
Let
$N_x:=N_{\partial A_{\varphi(x)},x},
$
where $A_{\varphi(x)} =\{y\colon\, \varphi(y) \le \varphi(x)\}$ and
$$
S:= \{x\in A\colon\,  N_x\ \hbox{contains more than one element}\},
$$
i.e. $S$ is the set of all the points $x$ such that the boundary of the
sub-level set containing $x$ is not differentiable at~$x$.
Then  $\mathcal{H}^{d}(S)=0$.
\end{lemma}

\begin{proof}
First we consider the case $d=2$.
Fix an orthonormal basis $\{e_1, e_2\}$
and identify every unit vector ${\rm n}$ with $\alpha \in [0,2\pi)$,
where $\alpha$ is the angle between $e_1$ and~${\rm n}$. We write
$$
{\rm n}:={\rm n}_{\alpha}: = \cos \alpha \cdot e_1 + \sin \alpha \cdot e_2.
$$
The set $S$ is a countable union of the sets
$$
S_{p,q} := \bigl\{ x\colon\ [p-q,p+q]\subset N_x\bigr\}, \quad
p,q \in \mathbb{Q}\cap [0,2\pi).
$$
If $S$ has a positive measure,
then  $\mathcal{H}^{d}(S_{p,q}) > 0$  for some $p,q$.
If $x \in S_{p,q}$, then we have $A_{\varphi(x)} \subset
\{z\colon\, \langle z-x,{\rm n}_{p}\rangle \le 0\}$.
Note that the line
$$
l_{x,p}(z) = \{z\colon\,
\langle z-x,{\rm n}_{p}\rangle = 0\}
$$
intersects $S_{p,q}$ exactly at $x$. Indeed, otherwise we get
two points $x,y$ such that the sub-level sets $A_{\varphi(x)}$ and
$A_{\varphi(y)}$ both intersect $l_{x,p}$ at two different points and belong
to the same half-plane $P$ with $\partial P = l_{x,p}$.
Hence neither $A_{\varphi(x)} \subset A_{\varphi(y)}$  nor
$A_{\varphi(y)} \subset A_{\varphi(x)}$ hold, which is impossible.
Thus we obtain that
$l_{x,p} \cap S_{p,q} = \{x\}$.
Finally, applying Fubini's theorem and disintegrating
Lebesgue measure along the lines parallel to $l_{0,p}$,
we obtain $\mathcal{H}^2 (S_{p,q})=0$, which is a contradiction.
So the lemma is proved for $d=2$.
The multi-dimensional case follows by induction and Fubini's theorem.
Indeed, fix an orthonormal basis $\{e_1, \ldots, e_d\}$.
Note that all sections of a convex body are convex. Disintegrating $S$ along
$e_i$ and applying the result for $d-1$, we obtain
that
$$S_i = \{\hbox{projection of $N_x$ on $x_i=0$ has more than one element}  \}
$$
has measure zero. Since $S= \bigcup_{i=1}^d S_i$, the proof is complete.
\end{proof}

{\bf Proof of Theorem \ref{maintheorem}:}
According to Lemma \ref{uniform-conv} and Lemma \ref{powergrad},
we have  $\varphi_n \to \varphi$ uniformly
and $|\nabla \varphi_{n}|^{\frac{1}{1+t_n}} \to 1$ almost everywhere.
It remains to prove that
$$\nabla \varphi_n/|\nabla \varphi_n| \to {\rm n}$$
almost everywhere.
Let us fix $x \in A$. Since $\varphi_n(x) \to \varphi(x)$, one has
$I_{A_{\varphi_n(x),n}} \to I_{A_{\varphi(x)}}$ almost everywhere, where
$$
A_{\varphi_n(x),n} = \{y\colon\, \varphi_n(y) \le \varphi_n(x)\},
\
A_{\varphi(x)} = \{y\colon\, \varphi(y) \le \varphi(x)\}.
$$
According to Lemma \ref{normalcone}, ${\rm n}(x)$ is well-defined for almost all $x$.
The same holds for every $\nabla \varphi_n(x)/|\nabla \varphi_n(x)|$.
So, without loss of generality we can fix $x$ in the interior
of $A$ such that ${\rm n}(x)$ and $\nabla \varphi_n(x)/|\nabla \varphi_n(x)|$
are well-defined.
If the vectors
$\nabla \varphi_n(x)/|\nabla \varphi_n(x)|$ do not converge to ${\rm n}(x)$,
then,
extracting a convergent
subsequence from a sequence of unit vectors
$\{\nabla \varphi_n(x)/|\nabla \varphi_n(x)|\}$,
we obtain a unit vector $\eta \ne {\rm n}(x)$.
By using convergence $I_{A_{\varphi_n(x),n}} \to I_{A_{\varphi(x)}}$,
one can show that $\langle\eta, z-x\rangle\le 0$ for all $z\in A$,
i.e., $\eta\in N_x$, which contradicts the choice of $x$.

It remains to verify the evolution equation for a smooth potential $\varphi$.
Indeed, let us choose an orthonormal basis $\{e_i\}$ at $x$ such that
$e_1={\rm n}$ and every vector $e_i$, $2 \le i \le d$, belongs to the
tangent space of $\partial A_{t}$
at $x$.
Let us write the change of variables formula for $T = \varphi \cdot {\rm n}$.
Differentiating along ${\rm n}$ we find
$$
\partial_{\rm n} T = \partial_{\rm n} \varphi \cdot {\rm n}
 + \varphi \cdot \partial_{\rm n} {\rm n}.
$$
Differentiating the identity $\langle {\rm n},{\rm n} \rangle  =1$,
we see that $\partial_{\rm n} {\rm n}$ belongs to the tangent space
of $\partial A_{t}$ at $x$. In addition,
$\partial_{\rm n} \varphi = |\nabla \varphi|$.
Next we note that
$$
\partial_{e_i} T =\varphi \cdot \partial_{e_i} {\rm n},
\quad
\langle \partial_{e_i}{\rm n},{\rm n}\rangle=0,
\quad
1\le i\le d.
$$
Hence
$$\det DT= |\nabla \varphi| \varphi^{d-1} \det
\bigl( \langle \partial_{e_i} n, e_j \rangle \bigr).
$$
Since
$
K = \det \bigl( \langle \partial_{e_i} {\rm n}, e_j \rangle,
$
we have
$\det DT=   |\nabla \varphi| \varphi^{d-1} K$.
Thus one obtains the following change of variables
formula (the Monge--Amp{\`e}re equation):
$$
\varrho_0 = \varrho_1( \varphi \cdot {\rm n}) |\nabla \varphi| \varphi^{d-1} K.
$$
It remains to note that the level sets $\partial A_s$ are shrinking with the velocity
$1/|\nabla \varphi|$ in the direction of $-{\rm n}$.
Hence (\ref{GaussFlow}) follows  from the change of variables formula.
The proof is complete.

\begin{example}
Let $A$ be a convex compact set.
Set
$$
\varrho_1(x) := \frac{C_{d,r}}{|x|^{d-1}} , \ \varrho_0(x) :=
\frac{1}{\mathcal{H}^d(A)},
$$
where ${\displaystyle C_{d,r} = \Bigl(\int_{B_r} \frac{dx}{|x|^{d-1}}\Bigr)^{-1}}$.
Varying $r$ we can show the existence of a weak solution
(in the ``transportation sense'') to the classical Gauss
curvature flow which starts from $\partial A$
and satisfies the equation
$$\dot{x}(s) = - c \ K(x)\cdot {\rm n}(x),$$
where $c$ can be chosen arbitrarily.

Certainly, a rigorous justification of this formula requires some
additional work, since we
have not proved that $\varphi$ is differentiable.
\end{example}

\section{Injectivity and uniqueness}

In this section, we prove that $T$ is invertible and essentially
unique.

Recall that the Legendre transform of a convex function $W$ on a
convex set $A$ is defined by
$$
W^{*}(y) = \sup_{x \in A} \bigl( \langle x,y \rangle - W(x)
\bigr).
$$
Let $\partial W(x)$ denote the subdifferential of $W$ at~$x$.
Recall also the following known fact from the
theory of convex functions (see, e.g., \cite[Theorem 23.5]{Rock}).

\begin{lemma}\label{subdiff}
Let $v \in \partial W(x)$ for every $x \in [x_1,
x_2]$, where $x_1 \ne x_2$ and $[x_1, x_2] = \{ t x_1 + (1-t) x_2,
\, t \in [0,1]\}$. Then $[x_1, x_2] \subset \partial W^{*}(v)$. In
particular, $W^{*}$ is not differentiable at~$v$.
\end{lemma}

In addition to the singular set
$S\subset A$ of all points $x$ such that
$N_{\partial A_t,x}$, where $t=\varphi(x)$, contains more than one
element, we introduce another set of degeneracy of $\rm{n}$ defined by
$$
U = \{ x \in A \backslash S\colon\, \hbox{there is
 $x' \in \partial A_{t}$, $t=\varphi(x)$,
 such that $x'\not=x$ and ${\rm{n}}(x) \in
N_{\partial A_t,x'}$} \}.
$$

\begin{proposition}\label{zerosets}
{\rm(i)}
 Consider the set $C= \partial A_t$ for some
fixed~$t$. Then
the set ${\rm{n}}(U \cap C)$
 in $S^{d-1}$ has  $\mathcal{H}^{d-1}$-measure zero.

{\rm(ii)}
The sets $T(U)$ and
$$
\widetilde{T}(S) := \bigcup_{x \in S} \varphi(x) \cdot N_{\partial
A_{\varphi(x)}, x}
$$
have $\nu$-measure zero.
\end{proposition}
\begin{proof}
(i)
It is sufficient to prove our claim locally on $C$ in a small
neighborhood $\mathcal{O}$ of a point $x_0$
where ${\rm{n}}(x_0)$ is unique. We may assume that
${\rm{n}}(x_0) = - e_d$,
 the surface $C\cap \mathcal{O}$ is the graph of a convex function
$W\colon\, B\subset \mathbb{R}^{d-1} \to \mathbb{R}$,
where $B$ is an open ball containing~$0$,
and that $W$ attains minimum at~$0$.
In addition, we may assume that
$\partial W(B)$ is a bounded set. We parameterize $C\cap \mathcal{O}$
in the following way:
$$
B\ni (x_1, \ldots, x_{d-1} ) \mapsto (x_1, \ldots, x_{d-1}, W(x)).
$$
Since $W$ is Lipschitzian on $B$, the surface measure
$\mathcal{H}^{d-1}$ on $C\cap\mathcal{O}$ corresponds to the measure
$(1+ |\nabla W|^2)^{\frac{1}{2}} \mathcal{H}^{d-1}$ on $\mathbb{R}^{d-1}$.
The Gauss map $\rm{n}$ is  given by
$$
{\rm{n}} = \frac{1}{\sqrt{1+|\nabla W|^2}} (-\partial_{x_1} W,
\ldots, -\partial_{x_{d-1}} W, 1 ).
$$
This holds for every $(x_1, \ldots, x_{d-1})\in B$ such that
$(x_1, \ldots, x_{d-1}, W(x)) \notin S \cap C\cap \mathcal{O}$.
The projection of
$S \cap C\cap \mathcal{O}$ on $B$ coincides with the points of
non-differentiability of~$W$.

It is convenient to identify the half-sphere $S^{d-1} \cap \{y_d
\le 0\}$ with its projection $\Pi^{d-1}$ on $\mathbb{R}^{d-1}$
 and $\rm{n}$ with the mapping $\widetilde{\rm{n}}\colon\,
-\frac{\nabla W}{\sqrt{1+ |\nabla W|^2}}$ taking values in
$\Pi^{d-1}$. Note that the surface measure on $S^{d-1}$ has the form
$m_{d-1}: = \frac{1}{\sqrt{1-|y|^2}}  \mathcal{H}^{d-1}$
in the local chart (on the set where $y_1^2+\cdots+y_{d-1}^2<1$)
$$
(y_1, \ldots, y_{d-1}) \mapsto \Bigl(y_1, \ldots, y_{d-1},
\sqrt{1-y_1^2 - \cdots  -y^2_{d-1}}\,\Bigr).
$$
Hence we have to show that
\begin{equation}\label{26.03.08}
m_{d-1}(F \circ \nabla W(U'))=0,
\end{equation}
where $U'$ is the corresponding projection of $U \cap C$ and
$$
F(x) =-\frac{x}{\sqrt{1+|x|^2}}.
$$
The mapping $F$ is smooth and nondegenerate everywhere. Hence
in order to prove (\ref{26.03.08}) it suffices to show that
$
\mathcal{H}^{d-1} (\nabla W(U'))=0.
$
Let us set $W:= +\infty$ outside of $B$. The Legendre transform
$W^{*}$ is finite everywhere. By Lemma \ref{subdiff}, the set $\nabla
W(U')$ is contained in the set of nondifferentiability of $W^{*}$,
hence has $\mathcal{H}^{d-1}$-measure zero.

(ii) By Fubini's theorem, it suffices to show that
for each $t$ the intersection of
the  set $T(U)$  with the
sphere of radius $t$ has zero $\mathcal{H}^{d-1}$-measure.
By construction, these intersection coincide with the sets
 $T(\partial A_t\cap U)$
  defined similarly.
 Therefore, the claim for $T(U)$ follows by assertion~(i).

In order to see that the set $\widetilde{T}(S)$ has $\nu$-measure zero,
we observe that its intersection with the  set
$T(A\backslash S)$ of full $\nu$-measure belongs to~$T(U)$,
which is clear from the definition of~$U$.
\end{proof}

Now we can show that $T$ is invertible.

\begin{corollary}
The mapping $T$ is injective on a set of full $\mu$-measure.
Hence there exists a measurable mapping
$T^{-1}\colon\, B_r \to A$ such that
$T(T^{-1}(y)) =y$ for
$\nu$-almost all $y$ and $T^{-1}(T(x)) = x$ for $\mu$-almost all~$x$.
\end{corollary}
\begin{proof}
Since the equality $T(x_1)=T(x_2)$ may only happen if
$\varphi(x_1)=\varphi(x_2)$, i.e., $x_1$ and $x_2$ belong to the
same level set $\partial A_t$,
it follows from our previous considerations that $T$ is injective
outside the set $T^{-1}(\widetilde{T}(S)\cup T(U))$. This set
has $\mu$-measure zero because the set $\widetilde{T}(S)\cup T(U)$
 has $\nu$-measure zero by the above proposition.
\end{proof}

\begin{theorem}
The mapping $T$ constructed above is unique in the following sense{\rm:}
if a measurable mapping $T_0\colon\, A\to B_r$ is such that
$\nu = \mu \circ T_0^{-1}$ and
$T_0 = \varphi_0 \cdot {\rm{n}}_0$, where
$\varphi_0\colon\, A\to [0,r]$ is a continuous function with
convex sub-level sets $A_{t,0}:=\{\varphi_0\le t\}$ and
${\rm n}_0$ is the corresponding
Gauss map, then $T=T_0$\ $\mu$-a.e.
\end{theorem}
\begin{proof}
Let us show that $\varphi_0(x)\le \varphi(x)$ for all $x\in A$. This will yield the equality $\varphi_0=\varphi$ because
otherwise there is $t$ such that
$\mu(\{\varphi_0\le t\})>\mu(\{\varphi\le t\})$, which is
impossible since both sides equal $\nu(B_t)$.
Set
$$
C_{t} := \Bigl\{x \in A\colon\ x \in
\partial A_{t,0} \cap A_{t}  \Bigr\},
\quad
D_{t} := \Bigl\{x \in A\colon\ x \in
\partial A_{t} \backslash {\rm{Int}}(A_{t,0}) \Bigr\},
$$
$$
U_{\tau} := \bigcup_{t \ge \tau} C_{t},
\quad
 V_{\tau} := \bigcup_{t \ge \tau} D_t.
$$
We observe that
for every $x \in \partial A_{t} \backslash {\rm{Int}} (A_{t,0})$
there exists $t' \ge t$ depending on $x$ such that
$x \in  C_{t'}$. Indeed,
if $x \in \partial A_{t,0} \cap \partial A_{t}$,
then $x \in C_t$. Otherwise one has
 $x \in \partial A_{t',0}$ for some $t'=t'(x) > t$.
 Since $x\in \partial A_{t}$, we have  $x\in {\rm{Int}}(A_{t'})$.
 Hence
\begin{equation}\label{v-sub-u}
V_\tau \subset U_\tau.
\end{equation}
For every Borel set $C \subset A$, set
$$
\widetilde{T}(C) := \bigcup_{x \in C}  \varphi(x) \cdot
N_{\partial A_{\varphi(x)},x}.
$$
Let us show that
\begin{equation}
\label{inclusion}
T_0(C_t)\subset \widetilde{T}(D_t ).
\end{equation}
Suppose that $x_0 \in \partial A_{t,0} \cap A_{t}$, $t=\varphi_0(x_0)$,
$v={\rm n}_{x_0,0} \in N_{\partial A_{t,0},x_0}$.
We show that $v \in T(\partial A_{t} \backslash {\rm{Int}} (A_{t,0}))$.
Let us consider
the support hyperplane $L_{x_0,v}\perp v$
to $A_{t,0}$ at~$x_0$. If $x_0 \in \partial A_{t}$
and $v \in N_{\partial A_t,x_0}$, the claim is obvious.
 Otherwise $L_{x_0,v}$ splits $A_{t}$ in two convex parts
$A'_{t}$ and $A^{''}_{t}$.
Since  $L_{x_0,v}$ is a support hyperplane to $A_{t,0}$,
one of these parts, say, $A^{''}_{t}$, and $A_{t,0}$
are separated by $L_{x_0,v}$.
There exists a hyperplane $L$  parallel to
$L_{x_0,v}$ that is supporting to $A^{''}_{t}$ and passes through a
point $x_1 \in \partial A^{''}_{t}$. Then
$v \in N_{\partial A_t,x_1}$. This proves (\ref{inclusion}).
Hence we have
\begin{equation}\label{u-v-v}
T_0(U_t) \subset \widetilde{T}(V_t),\quad 0 \le t \le r.
\end{equation}
Suppose now that there exists $x_0$ such that
$\varphi_0(x_0)>\varphi(x_0)$. Then,
by the continuity of $\varphi$ and $\varphi_0$,  there is
$\tau>0$ for which the inclusion in
(\ref{v-sub-u}) is strict and there is a neighborhood in $U_\tau$ not intersecting $V_\tau$. Therefore,
$\mu(U_{\tau})>\mu(V_{\tau})$.
Taking into account that $T$ is injective on a full measure set,
we obtain
$$
\nu(T_0(U_{\tau}))=
\mu(T_0^{-1}(T_0(U_{\tau})))
 \ge \mu(U_{\tau}) > \mu(V_{\tau}) = \nu(T(V_{\tau})),
$$
which contradicts (\ref{u-v-v}) because
$ \nu(T(V_{\tau}))= \nu(\widetilde{T}(V_{\tau}))$
according to Corollary \ref{zerosets}.
\end{proof}

\section{Duality}

Now we  consider
certain duality properties of the potential~$\varphi$.
The duality principle of Kantorovich is a powerful tool for investigating
the Monge--Kantorovich problem.
In our case we also have a kind of the duality formula which
relates the potential
$\varphi$ to  some function $\psi$ that can be considered as the
support function
of the family of level sets~$A_t$.
Note that some interesting duality results
for the solution of the Monge--Kantorovich problem on a
sphere  with applications
to the prescribed Gauss curvature problem
have been obtained in~\cite{Oliker}.

For every $y \in B_r$ we set

\begin{equation}
\label{support_f}
\psi(y) = \sup_{x\colon\, \varphi(x) \le |y|} \langle x, y \rangle.
\end{equation}
Note that the restriction of $\psi$ to
$\partial B_{|y|}$ coincides with the support function $S_{A_{|y|}}$ of
$A_{|y|} = \{x\colon\, \varphi(x) \le |y|\}$, where
the support function is defined by
$$
S_{A_{|y|}}(v) := \sup_{x \in A_{|y|}} \langle v, x \rangle.
$$

\begin{lemma}
For $\nu$-almost all $y$ one has
\begin{equation}
\label{ST}
\psi(y) = \langle T^{-1}(y), y\rangle.
\end{equation}
\end{lemma}
\begin{proof}
It is clear that the supremum on the right-hand side
of (\ref{support_f}) is attained at a point $p$
such that $y \in N_{\partial A_{\varphi(p)},p}$.
This implies that $p$ coincides with $T^{-1}(y)$
 for $\nu$-almost all~$y$, hence is $\nu$-almost everywhere
  well-defined, which yields our claim.
\end{proof}

Now we show how to describe $\psi$ as a limit of certain functions
depending on pre-limit potentials~$\varphi_t$.
Recall that the Legendre transform $W^{*}_t$ satisfies
the inequality
\begin{equation}
\label{Legendre}
W_t(x) + W^{*}_t(y) \ge \langle x, y \rangle.
\end{equation}
An equality holds if and only if
$y \in \partial W_t(x)$ and $x \in \partial W^*_t(y)$.
Moreover, $W_t$ and $W^{*}_t$
satisfy the identities
$$
\nabla W^*_t \circ \nabla W_t (x)  =x, \quad
\nabla W_t \circ \nabla W^*_t (y)  =y
$$
almost everywhere on the sets $A$ and $\nabla W_t(A)$.
Since
$$
\nabla W_t = |T_t|^t T_t,
$$
one has
$$
T^{-1}_t(y) = \nabla W^*_t (|y|^t y).
$$
In what follows we denote by $I$ the identity matrix and by $I_z$  the
orthogonal projector on the one-dimensional vector subspace
generated by $z$, i.e.,
$$
I_z v = \frac{\langle v, z\rangle}{|z|} \frac{z}{|z|}.
$$

We have found a sequence $t_n\to +\infty$ for which the mappings $T_{t_n}$
converge to $T$ almost everywhere on~$A$, hence converges
in measure~$\mu$. For this sequence, the
following holds.

\begin{lemma}
\label{inv-conv}
The mappings $T_{t_n}^{-1}$ converge to $T^{-1}$ in measure~$\nu$.
Hence there exists a subsequence $t_n' \to \infty$ such that
$T^{-1}_{t_n'} \to T^{-1}$
\ $\nu$-almost everywhere.
\end{lemma}
\begin{proof}
Since $\mu\circ T^{-1}=\mu\circ T_{t_n}^{-1}=\nu$,
for any $\nu$-measurable function~$f$, the functions
$f\circ T_{t_n}$ converge to $f\circ T$ in measure~$\mu$
(see, e.g., \cite[Corollary~9.9.11]{B}).
Hence the mappings $T^{-1}\circ T_{t_n}$ converge
to $T^{-1}\circ T=I$ in measure~$\mu$. Therefore,
for every $c>0$ one has
$$
\nu\Bigl(y\colon\ |T_{t_n}^{-1}(y)-T^{-1}(y)|\ge c\Bigr)
=
\mu\Bigl(x\in A\colon\ |x-T^{-1}T_{t_n}(x)|\ge c\Bigr)\to 0
$$
as $n\to\infty$, which completes the proof.
\end{proof}

\begin{theorem}
Let a function $\psi_t$ be defined by the relation
$$
W^*_t(z) = |z|^{\frac{t}{1+t}} \psi_t \bigl(z |z|^{-\frac{t}{1+t}} \bigr).
$$
Equivalently,
$$
\psi_t(y) = \frac{W^{*}_t(y|y|^t)}{|y|^t}.
$$
Then one has $\psi = \lim\limits_{t_n \to \infty} \psi_{t_n}$
almost everywhere for some sequence $\{t_n\}$.
\end{theorem}
\begin{proof}
Note that it is consistent with our previous choice
of $W_t$ to assume that $W^{*}_t(0)=0$. Indeed,
$W_t(x) \ge \langle x, y \rangle - W^{*}_t(y)$, hence
taking $y=0$ we find $W_t(x) \ge 0$.
Taking any $x_0 \in \partial W^{*}_t(0)$ we easily obtain
 $W_t(x_0)=0$. Indeed, for $(x_0, 0)$
inequality (\ref{Legendre}) becomes an equality, hence
$
W_t(x_0) + W^{*}_t(0) = \langle x_0, 0 \rangle =0.
$

The inequality
$W^{*}_t(a) - W^{*}_t(b) \le \langle a-b, \nabla W^{*}(a)\rangle$
yields, by substituting $b=0$ and $a=y|y|^t$, that
$$
\psi_t(y) \le \langle \nabla W^{*}_t(y|y|^t), y \rangle = \langle y, T^{-1}_t(y) \rangle.
$$
Similarly, if $a=0$ and $b=y|y|^t$, one has
 $\psi_t(y) \ge \langle v, y\rangle$ for
 any $v \in \partial W^{*}_t(0)$. In particular,
$$
|\psi_t(y)| \le {\rm{diam}}(A) |y|.
$$
One has
$$
\nabla W^*_t (z)
= {\frac{t }{1+t}}|z|^{-\frac{1}{1+t}}
\frac{z}{|z|} \cdot \psi_t \bigl(z |z|^{-\frac{t}{1+t}} \bigr)
+
\Bigl(I - \frac{t}{1+t}I_z\Bigr)
\nabla \psi_t \bigl(z |z|^{-\frac{t}{1+t}} \bigr).
$$
Substituting $z = |y|^t y$ we obtain
\begin{equation}
\label{TP}
T^{-1}_t(y)
={\frac{t }{1+t}}  \psi_t(y) \frac{y}{|y|^2}
+
\Bigl(I - \frac{t}{1+t}I_y\Bigr) \nabla \psi_t (y).
\end{equation}
Taking the inner product with $y$ we find
\begin{equation}
\label{TtSt}
\langle T^{-1}_t(y), y \rangle
={\frac{t }{1+t}}  \psi_t (y)
+\frac{1}{1+t} \langle \nabla \psi_t (y), y \rangle.
\end{equation}
In view of Lemma \ref{inv-conv} and
equality (\ref{ST}) it suffices to show that
$\frac{1}{1+t_n}\langle \nabla \psi_{t_n} (y), y \rangle \to 0$ almost everywhere for some $\{t_n\}$.
Indeed, since $\psi_t \le \langle T^{-1}_t(y), y\rangle$,
 we obtain from (\ref{TtSt}) that
$$
\langle T^{-1}_t(y), y \rangle
\le
{\frac{t }{1+t}}  \langle T^{-1}_t(y), y\rangle
+
\frac{1}{1+t} \langle \nabla \psi_t (y), y \rangle.
$$
Hence
$$
\psi_t \le \langle T^{-1}_t(y), y\rangle \le
\langle \nabla \psi_t (y), y \rangle.
$$
Taking into account that $\psi_t \ge -{\rm{diam}}(A) |y|$,
 we see that $\{\langle \nabla \psi_t (y), y \rangle\}$
is uniformly bounded from below.
The integration by parts formula yields
$$
\int_{B_r} \langle \nabla \psi_t (y), y \rangle\, dx
= -d \int_{B_r} \psi_t\, dx
+ \int_{\partial B_r} |y| \psi_t\, d\mathcal{H}^{d-1}.
$$
Applying again the estimate $|\psi_t| \le {\rm{diam}}(A) |y|$ we obtain
$$
\sup_t \int_{B_r} \langle \nabla \psi_t (y), y \rangle\, dx < \infty,
$$
hence
$\sup_{t} \| \langle \nabla \psi_t (y), y \rangle \|_{L^1(B_r)} < \infty$.
Therefore,
$$
\lim\limits_{t \to \infty} \frac{1}{1+t}\|
\langle \nabla \psi_t (y), y \rangle \|_{L^1(B_r)} =0.
$$
 Extracting a subsequence we complete the proof.
\end{proof}

\begin{remark}
Taking a scalar product of (\ref{TP}) with
any vector $v \bot y$ we obtain the equality
$\partial_v \psi_n(y) =
 \langle T^{-1}_n(y), v \rangle$.
 Let us set
 $$
 \partial_v \psi(y) := \lim_{t_n \to \infty} \partial_v \psi_n(y).
 $$
 In view of convergence $T_n \to T$ this definition makes sense.
 Moreover, we  have
 \begin{equation}
\label{S'}
\partial_v \psi(y) =
 \langle T^{-1}(y), v \rangle,
 \end{equation}
for any $v \bot y$.

Taking into account (\ref{ST}) we obtain the following remarkable relation:
$$
T^{-1}(y)
= \frac{\psi(y)}{|y|} e_1(y)
+
\sum_{i=2}^{d} \partial_{e_i(y)} \psi(y) \ e_i(y),
$$
where $\{e_i(y)\}$ is an orthonormal system of unit
vectors chosen in such a way that
 $e_1(y) = y/|y|$ and $e_i(y) \bot y$, $2 \le i \le d$.
\end{remark}

\begin{remark}
Let us see what happens in the limit with  the duality formula
$$
W_t(x) + W^{*}_t(z) \ge \langle x, z \rangle.
$$
It can be rewritten as
$$
\frac{1}{t+2} \varphi^{t+2}_t(x) + |y|^{t} \psi_t(y)
\ge \langle x,y \rangle |y|^t
$$
by letting $z:=y|y|^t$.
If
$y = T_t(x) = \varphi_t
\frac{\nabla \varphi_t}{|\nabla \varphi_t|}
|\nabla \varphi_t|^{\frac{1}{1+t}}$,
then an equality holds. It is known that this is possible
only if the pair $(x,y)$ belongs to the graph of~$T_t$.
Hence
we obtain the following duality relation:
$$
\frac{1}{t+2} \varphi^2_t(x)
|\nabla \varphi_t(x)|^{-\frac{t}{t+1}} +
\psi_t(T_t(x)) = \langle x, T_t(x) \rangle.
$$
In the limit $t \to \infty$ we find
$$
\psi(T(x)) = \langle x, T(x) \rangle.
$$
\end{remark}

It is worth noting that for the constructed transformation $T$
one can establish a change of variables formula involving
certain analogs of Alexandroff's determinants.
Such formulas which neglect singular components are known
for optimal transformations and triangular transformations
(concerning the latter, see \cite[Ch.~10]{B} and \cite{BKM}).

\vskip .2in

This work was supported by the RFBR projects 07-01-00536,
08-01-91205-JF, GFEN-06-01-39003,
RF Prezident Grant MD-764.2008.1,
DFG Grant 436 RUS 113/343/0(R), ARC Discovery Grant DP0663153, and
the SFB 701 at the University of Bielefeld.
The second-named author thanks Franck Barthe for his hospitality and fruitful
discussions during the author's visit to the University of Paul Sabatier in
Toulouse, where this work was partially done.

\end{document}